\documentclass{amsart}
\usepackage{amsmath,amsfonts,amssymb,amsthm}
\usepackage{hyperref}

\newcommand{\Ff}[1][l]{\mathbf{F}_{#1}}
\newcommand{\GL}{\mathrm{GL}}
\newcommand{\Q}{\mathbf{Q}}

\theoremstyle{plain}

\newtheorem*{conj}{Conjecture}

\theoremstyle{remark}

\theoremstyle{definition}

\title{On some remarkable congruences between two elliptic curves}
\author{Nicolas Billerey}
\address{(1) Universit\'e Clermont Auvergne, Universit\'e Blaise Pascal,
	Laboratoire de Math\'ematiques,
	BP 10448,
	F-63000 Clermont-Ferrand, France.
	(2) CNRS, UMR 6620, LM, F-63171 Aubi\`ere, France}
\email{Nicolas.Billerey@math.univ-bpclermont.fr}
\date{}

\begin{document}

\begin{abstract}
We exhibit two non-isogenous rational elliptic curves with \mbox{$17$-torsion} subgroups isomorphic as Galois modules.
\end{abstract}

\maketitle

\section*{}
Consider the following rational elliptic curves 
\[
E\colon y^2+xy=x^3-8x+27
\]
and 
\[
E'\colon y^2+xy=x^3+8124402x-11887136703
\] 
labelled \cite[\href{http://www.lmfdb.org/EllipticCurve/Q/3675.g1}{3675.g1}]{lmfdb} and~\cite[\href{http://www.lmfdb.org/EllipticCurve/Q/47775.be1}{47775.be1}]{lmfdb} in LMFDB res\-pectively. Their discriminants are given by the following identities~:
\[
\Delta(E)=-3^5\cdot 5^2\cdot 7^2\quad\text{and}\quad \Delta(E')=-3^5\cdot 5^2\cdot 7^2\cdot 13^{17}.
\]

The aim of this note is to prove that their \(17\)-torsion groups are isomorphic as Galois modules\footnote{As pointed out to the author, these congruences  were known to John Cremona since 2007 (unpublished).}. We apply Proposition~4 of~\cite{KrOe92} with (in their notation) \(N=3675\) and \(N'=47775\). We then have~\(M=N'\) (both~\(E\) and~\(E'\) have bad multiplicative split reduction at~\(3\)) and \(\mu(M)=20160\). We thus have to check that
\[
a_l(E)\equiv a_l(E')\pmod{17}\quad\text{for every prime \(l<3360\), \(l\not= 3,5,7,13\)}
\]
and
\[
a_{13}(E)a_{13}(F)\equiv\pm(13+1)\pmod{17}.
\]
This can be easily done using~\textsc{Sage} (\cite{sage}). Therefore, the mod 17 representations attached to~\(E\) and~\(E'\), denoted (as in~\cite{KrOe92}) \(\rho\) and~\(\rho'\) respectively, have isomorphic semi-simplifications. However, both representations are irreducible. This follows from the fact that both curves have bad additive reduction at~\(5\) with \(\mathrm{ord}_5(\Delta(E))=\mathrm{ord}_5(\Delta(E'))=2\).  Hence the group \(\Phi_5\) is cyclic of order~\(6\) and the images of~\(\rho\) and~\(\rho'\) cannot be included in a Borel subgroup of~\(\GL_2(\Ff[17])\). For more details and the definition of the group~\(\Phi_5\), see \S5.6 and Ex.~5.7.4 of~\cite{Ser72} (or Prop.~3.3 of~\cite{Bil11}). 

As of April 2016, there are eight known pairs (up to isomorphism) of non-isogenous  $17$-congruent rational elliptic curves with conductors less than $360,000$. They turn out to be all quadratic twists of the pair~$(E,E')$. Besides, to the author's knowledge there is no known example of two non-isogenous rational curves with $p$-torsion subgroups isomorphic as Galois modules for some prime number~$p>17$. Any such example would be quite interesting in view of the conjecture below.
\begin{conj}[Frey-Mazur]
There exists a constant~\(C\) such that for any prime $p\geq C$ and any pair of elliptic curves~\(E,E'\) over~\(\Q\), the following holds~:
\[
E[p]\text{ and }E'[p]\text{ are isomorphic as Galois modules }\Longrightarrow E\text{ and \(E'\) are isogenous}.
\]
\end{conj}

\newcommand{\etalchar}[1]{$^{#1}$}

%
\end{document}